\documentclass[12pt]{article}
\usepackage{fancyhdr, array,calc,graphicx,url,tabularx}
\usepackage{geometry}
\usepackage{latexsym}
\usepackage{amssymb}
\usepackage{amsmath}
\usepackage{makeidx}

\makeindex

{
   \end{minipage}
   \vspace*{\stretch{3}}
   \clearpage
}

\newcommand{\bR}{{\mathbb{R}}}

\newcommand{\rest}{\restriction}

\newcommand{\card}[1]{{\vert #1 \vert} }

\renewcommand{\models}{\vDash}
\newcommand{\powerset}{{\wp}}
\newcommand{\cp}{{\rm cp }}

\newtheorem{theorem}{Theorem}[section]

\newtheorem{definition}[theorem]{Definition}
\newtheorem{rem}[theorem]{Remark}

\newtheorem{lemma}[theorem]{Lemma}
\newtheorem{corollary}[theorem]{Corollary}

\newtheorem{conjecture}[theorem]{Conjecture}
\newtheorem{question}[theorem]{Question}

\numberwithin{figure}{section}

\newenvironment{proof}{{\it{
Proof.}}}{\nopagebreak\mbox{}{\hfill$\square$}
\par\bigskip}
\newenvironment{Claim}{\noindent{\it
Claim} {}}{}
\newenvironment{Case}{\noindent{\it
Case}}{}

\newcommand{\rthm}[1]{Theorem~\ref{#1}}
\newcommand{\rlem}[1]{Lemma~\ref{#1}}

\newcommand{\rcor}[1]{Corollary~\ref{#1}}

\newcommand{\rsubsec}[1]{Subsection~\ref{#1}}

\newcommand{\rrem}[1]{Remark~\ref{#1}}

\def\inseg{\trianglelefteq}

\def\k{\kappa}
\def\a{\alpha}
\def\b{\beta}
\def\d{\delta}

\def\l{\lambda}

\def\P{{\mathcal{P} }}

\def\Q{{\mathcal{ Q}}}

\def\R{{\mathcal R}}

\def\H{{\rm{HOD}}}
\def\M{{\mathcal{M}}}
\def\N{{\mathcal{N}}}

\def\U{{\mathcal{U}}}
\def\S{{\mathcal{S}}}

\def\VT{{\vec{\mathcal{T}}}}
\def\VU{{\vec{\mathcal{U}}}}

\def\cp #1{{ crit  #1 }}
\def\card#1{\left|#1\right|}

\def\iff{\mathrel{\leftrightarrow}}

\def\and{\mathrel{\kern1pt\&\kern1pt}}

\def\inseg{\triangleleft}
\def\insegeq{\trianglelefteq}

\def\<#1>{\langle\,#1\,\rangle}

 \input xy
 \xyoption{all}
 
\begin{document}

\title{Hod up to $\text{AD}_{\mathbb{R}}+\Theta \text{ is measurable}$
 }
\date{}
\author{Rachid Atmai\\
Department of Mathematics\\
University of North Texas\\
General Academics Building 435\\
1155 Union Circle $\#$311430\\
Denton, TX 76203-5017\\
atmai.rachid@gmail.com\\\\
Grigor Sargsyan\thanks{First author's work is partially based upon work supported by the National Science Foundation under Grant No DMS-1352034 and DMS-1201348.}\\
        Department of Mathematics\\
        Rutgers University\\
        110 Frelinghuysen Rd.\\
        Piscataway, NJ 08854\\
        http://math.rutgers.edu/$\sim$gs481\\
        grigor@math.rutgers.edu,\\\\
}

\maketitle
\begin{abstract}
Suppose $M$ is a transitive class size model of $AD_{\bR}+``\Theta$ is regular". $M$ is a minimal model of $AD_{\bR}+``\Theta$ is measurable" if (i) $\bR, Ord\subseteq M$ (ii) there is $\mu\in M$ such that $M\models ``\mu$ is a normal $\bR$-complete measure on $\Theta$" and (iii) for any transitive class size $N\subsetneq M$ such that $\bR\subseteq N$, $N\models ``$there is no $\bR$-complete measure on $\Theta$". Continuing Trang's work in \cite{TrangHod}, we compute $\H$ of a minimal model of $AD_{\bR}+``\Theta$ is measurable".
\end{abstract} 

The computation of $\H$ of models of determinacy has been one of the central themes in descriptive inner model theory. Steel's seminal \cite{Steel1995} jumpstarted the project and a later work of Steel and Woodin (for instance see \cite{ATHM} or \cite{BSL}) established connections with the Mouse Set Conjecture, core model induction and the inner model problem. 

The analysis of $\H$ presented in the above papers, however, only computes $V_{\Theta}^\H$. Woodin computed the full $\H$ of $L(\bR)$ under $AD^{L(\mathbb{R})}$ (the proof can be found in \cite{OIMT}). Trang continued this work in \cite{TrangHod}. He presented the exact inner model theoretic structure of $\H$ of models of determinacy that have the form $L(\powerset(\bR))$ and are contained inside the minimal model of $AD_{\mathbb{R}}+``\Theta$ is regular". We extend Trang's work to minimal models of $AD_{\bR}+``\Theta$ is measurable". This notion, however, needs a definition.

Below we say that a measure $\mu$ is $\bR$-complete if whenever $(A_x: x\in \bR)\subseteq \mu$, $\cap_{x\in \bR}A_x\in \mu$. 

\begin{definition}
Assume $AD^+$.  We say $``\Theta$ is measurable" if there is an $\mathbb{R}$-complete normal measure $\mu$ on $\Theta$. We then let $\Theta ms$ be the theory $AD_{\bR}+``\Theta$ is measurable".
\end{definition}

Because we need to assume that the computation of $V_\Theta^{\H^{L(\powerset(\bR))}}$ can be carried out, our current techniques can only work in some minimal setting, in a setting where the Mouse Set Conjecture (MSC, see Chapter 3.1 of \cite{BSL}) and The Generation Of Pointclasses (see Conjecture 3.17 of \cite{BSL}) are true. The first says that ordinal definable reals appear in mice and the second says that sufficiently closed pointclasses are generated by hod pairs. The models where our computation works are the minimal models of $\Theta ms$.

\begin{definition}\label{minimal models}
$M$ is a minimal model of $\Theta ms$ if $M\models \Theta ms$, $\bR\cup Ord\subseteq M$ and for any $N \subsetneq M$ such that $\mathbb{R}\cup Ord\subseteq N$, $N \models \neg \Theta ms$. 
\end{definition}

 We now have the following simple representation of minimal models of $\Theta ms$.

\begin{lemma}\label{characterization of min models}
Suppose that $V$ is a minimal model of $\Theta ms$. Let $\mu$ be a normal $\mathbb{R}$-complete measure on $\Theta$. Then $V=L(\powerset(\mathbb{R}))[\mu]$.
\end{lemma}
\begin{proof}
We have that $L(\powerset(\mathbb{R}))[\mu]\models \Theta ms$. Because $L(\powerset(\mathbb{R}))[\mu]\subseteq V$, we have that $V=L(\powerset(\mathbb{R}))[\mu]$. 
\end{proof}

It is not immediately clear that the existence of a minimal model of $\Theta ms$ follows from the existence of a model of $\Theta ms$. This is because there could be $\subsetneq$-descending sequence of models $(M_i: i<\omega)$ such that  for all $i<\omega$, $\Theta^{M_i}=\Theta^{M_{i+1}}$ and $M_i\models \Theta ms$. We will in fact show that this cannot happen (see \rthm{existence of minimal models}). However, towards showing the aforementioned fact, we will need to analyze the \textit{weakly minimal} models of $\Theta ms$. 
\begin{definition}
We say $M$ is a weakly minimal model of $\Theta ms$ if 
\begin{enumerate}
\item $\bR\cup Ord\subseteq M$ and $M\models \Theta ms$, 
\item for some $\mu\in M$ such that $M\models ``\mu$ is a normal $\bR$-complete measure on $\Theta$", $M=L(\powerset(\bR))[\mu]$, and
\item for any transitive $N$ such that $\powerset(\bR)^N\subsetneq \powerset(\bR)^M$ and $\mathbb{R}\cup Ord\subseteq N$, $N \models \neg \Theta ms$. 
\end{enumerate}
\end{definition}

It is important to note that it follows from the main results of \cite{ATHM} (see Section 6.4 and Section 6.1 of \cite{ATHM}) that if $M$ is a weakly minimal model of $\Theta ms$ then both MSC and The Generation of Pointclasses hold in $L(\powerset(\bR)^M)$. We will exploit this in the next section. We mention that while it is not clear that there are minimal models of $\Theta ms$, it is not hard to show, modulo the existing literature, that weakly minimal models do exist.

\begin{lemma}\label{weakly min models} Suppose $M$ is a transitive model of $\Theta ms$ such that $\bR, Ord\subseteq M$. Then there is a weakly minimal model of $\Theta ms$. 
\end{lemma}
\begin{proof}
Let $\Gamma\subseteq \powerset(\bR)^M$ be Wadge least initial segment of $\powerset(\bR)^M$ such that for some filter $\mu$ on $w(\Gamma)$, $L(\Gamma)[\mu]\models \Theta ms$\footnote{Recall that $w(\Gamma)$ is the supremum of the Wadge ranks of sets in $\Gamma$.}. Clearly, $L(\Gamma)[\mu]$ is a weakly minimal model of $\Theta ms$.
\end{proof}

\textbf{Acknowledgments.} The authors would like to thank Nam Trang for communicating them the problem considered in this paper. Also, the authors would like to thank the referee for invaluable comments.

\section{The main theorem}

Suppose $V$ is a weakly minimal model of $\Theta ms$. As was mentioned above, both MSC and The Generation of Pointclasses hold in $V$. We can then freely use the machinery of \cite{ATHM}. Let $M=L(\powerset(\bR))$ and working in $M$, let $\mathcal{F}$ be the set of hod pairs $(\mathcal{P},\Sigma)$ such that $\Sigma$ has branch condensation and is fullness preserving.

Following \cite{ATHM}, for $(\mathcal{P},\Sigma),(\mathcal{Q},\Lambda) \in \mathcal{F}$, we let \[(\mathcal{P},\Sigma) \preceq (\mathcal{Q},\Lambda)\]
if  for some $\a\leq \l^\Q$, $(\Q(\a), \Lambda_{\Q(\a)})$ is a tail of $(\P, \Sigma)$, i.e., $\Q(\a)\in I(\P, \Sigma)$ and $\Lambda_{\Q(\a)}=\Sigma_{\Q(\a)}$. It follows from comparison theory of hod pairs (see Chapter 2 of \cite{ATHM}) that $\preceq$ is directed. We then let $\mathcal{M}_{\infty}$ be the direct limit of the system $\mathcal{F}$ under the maps $\pi^{\Sigma}_{\mathcal{P},\mathcal{Q}}: \mathcal{P} \to \mathcal{Q}(\alpha)$ where $\a\leq\l^\Q$ is such that $\Q(\a)$ is a $\Sigma$-iterate of $\P$.

It is shown in \cite{ATHM} (see Theorem 4.24) that in $M$, $\M_\infty=V_\Theta^\H$. The following is essentially the generic interpretability result of \cite{ATHM} (see Theorem 3.10 of \cite{ATHM}). Given a hod premouse $\P$ and $X$ generic over $\P$, we let $\Sigma^{\P[X]}$ be the interpretation of $\Sigma^\P$ onto $\P[X]$ according to the procedure described in the proof of Theorem 3.10 of \cite{ATHM}.
 
\begin{lemma}\label{strategy for m_infty}
In $M$ (and hence in $V$), $\M_\infty$ is $(\Theta, \Theta)$-iterable via a strategy $\Sigma$ such that given any $\VT$ according to $\Sigma$,
\begin{center}
$\Sigma(\VT)=b \iff \M_\infty[\VT]\models \Sigma^{\M_\infty[\VT]}(\VT)=b$.
\end{center}
\end{lemma}

We let $\Sigma$ be the strategy of $\M_\infty$ described in \rlem{strategy for m_infty}. Next, we define a model extending $\M_\infty$. Given a hybrid premouse $\mathcal{N}\trianglerighteq \mathcal{M}_\infty$ we say $\N$ is \textit{good} if
\begin{enumerate}
\item  $\mathcal{N}$ is sound, and 
\item whenever 
\begin{center} 
$\pi: \bar{\mathcal{N}} \to \mathcal{N}$ 
\end{center}
is elementary and $\bar{\mathcal{N}}$ is countable, then $\bar{\mathcal{N}}$ is $\omega_1$-iterable above $\pi^{-1}(\Theta)$ as a $\Sigma^{\pi}$-premouse. 
\end{enumerate}

The next two lemmas are basic lemmas about good mice. 

\begin{lemma}\label{good mice line up} Suppose $\N_0, \N_1$ are good such that for some $\eta\in \N_0\cap \N_1$, $\eta$ is a cutpoint of both $\N_0$ and $\N_1$, $\N_0|\eta=\N_1|\eta$ and $\rho_\omega(\N_0), \rho_\omega(\N_1)\leq \eta$.  Then either $\mathcal{N}_0 \unlhd \mathcal{N}_1$ or $\mathcal{N}_1 \unlhd \mathcal{N}_0$. 
\end{lemma}
\begin{proof}
Suppose that neither $\mathcal{N}_0 \unlhd \mathcal{N}_1$ nor $\mathcal{N}_1 \unlhd \mathcal{N}_0$ holds. Let then $\pi: H \to L_{\xi}(\mathcal{P}(\mathbb{R}))[\mu]$ be elementary such that $\xi > \Theta$, $\mathcal{N}_i \in rng(\pi)$, $\Sigma \in rng(\pi)$, $\eta\in rng(\pi)$ and $\vert H \vert=\omega$.

 We let $\bar{\mathcal{N}}_i=\pi^{-1}(\mathcal{N}_i)$ for $i=0,1$. Because $\Sigma$ has hull condensation, it follows that $\pi^{-1}(\Sigma) =\Sigma^{\pi}\restriction H$. It follows from elementarity that for $i=0,1$, $\rho_{\omega}(\bar{\mathcal{N}}_i)=\pi^{-1}(\eta)$. But now because $\bar{\mathcal{N}}_i$ are sound $\Sigma^{\pi}$-mice, we have that $\bar{\N}_0\insegeq \bar{\N}_1$ or $\bar{\N}_1\insegeq \bar{\N}_0$.
\end{proof}

\begin{lemma} Suppose $\N$ is good such that $\rho_\omega(\N)\leq \Theta$. Then $\rho_\omega(\N)=\Theta$.
\end{lemma}
\begin{proof}
Towards a contradiction assume $\rho(\N)<\Theta$. Let $\pi: H\to L_{\xi}(\mathcal{P}(\mathbb{R}))[\mu]$ be such that $H\cap \Theta\in \Theta$, $\cp(\pi)>\rho(\N)$ and $\N\in rng(\pi)$. Let $\P=\pi^{-1}(\N)$. It follows from the proof of \rlem{good mice line up} that $M\models ``\P$ is OD". It then follows that $\P\insegeq \M_\infty$. This is a contradiction. 
\end{proof}

We then let $\M^*=\cup\{\N: \M_\infty\insegeq \N$, $\N$ is good and $\rho_\omega(\N)=\Theta\}$. Notice that we have that $\M^*\subseteq \H$. Let $\eta=o(\M^*)$.

Let now $\mu\in V$ be a normal $\bR$-complete measure on $\Theta$. Notice that if $V$ was a minimal rather than a weakly minimal model of $\Theta$ms then \rlem{characterization of min models} would imply that $V=L(\powerset(\bR))[\mu]$. Working in $\H_\mu$, let $\pi_\mu$ be the ultrapower embedding via $\mu\cap \H_\mu$. Let $\l_\mu=(\Theta^{++})^{\pi_\mu(\M^*)}$. Notice that the ordinal $\l_\mu$ may depend on $\mu$. Let then 
\begin{center}
$\M_\mu=\pi_\mu(\M^*)|\l_\mu$ 
\end{center}
and let $E_\mu$ be $(\Theta, \eta)$-extender derived from $\pi_\mu\rest \M^*$. 
 More precisely, 
\begin{center}
$(a, A)\in E_\mu \iff a\in (o(\M^*))^{<\omega}, A\in [\Theta]^{\card{a}}\cap \M^*,\ \text{and}\ a\in \pi_\mu(A)$
\end{center}
%

The following is our main theorem. 

\begin{theorem}[Main Theorem]\label{main theorem} Assume $V$ is a minimal model of $\Theta ms$. Suppose $\mu$ is a normal $\bR$-complete measure on $\Theta$. Then $\H=L[\M_\mu, E_\mu]$. 
\end{theorem}

We will present the proof as a sequence of lemmas. In the subsections that follows, we assume that $V$ is a minimal model of $\Theta ms$. However, we remark that

\begin{rem}\label{remark} all the results that follow except the results of \rsubsec{uniqueness of measures sec} can be carried out under the assumption that $V$ is just a weakly minimal model of $\Theta ms$. We ask the reader to keep this remark in mind while reading the subsequent sections.
\end{rem}

Before we go into the proof of the main theorem, we list some of the complications involved with proving it. First we will show that $E_\mu$ is amenable to $\M_\mu$ (see \rlem{amenability}). It follows from its definition that it coheres $\M_\mu$. It then follows that $(\M_\mu, E_\mu, \in)$ is a hod premouse. 

The next challenge is to show that no level of $L[\M_\mu, E_\mu]$ projects to or below $\Theta$ (see \rlem{not projecting}). A consequence of this is that $V_\Theta^{\H^M}=V_\Theta^{L[\M_\mu, E_\mu]}$. This then allows us to show that $\powerset(\bR)$ can be symmetrically added to $L[\M_\mu, E_\mu]$ (which is done as part of proving \rlem{not projecting}). It then follows that the model $L[\M_\mu, E_\mu]$ is independent of $\mu$. The final piece of the argument is that there is a unique normal $\bR$-complete measure $\mu$ on $\Theta$ minimizing $\lambda_\mu$. It then immediately follows that $L[\M_\mu, E_\mu]\subseteq \H$. Combining the aforementioned results it is then not hard to see that $V$ is a symmetric extension of $L[\M_\mu, E_\mu]$, which then easily implies that in fact $\H\subseteq L[\M_\mu, E_\mu]$.

\subsection{Amenability} 

Recall that $\eta=o(\M^*)$. We start by showing that 
\begin{lemma}\label{lower part characterization of m} $\M^*=\pi_\mu(\M_\infty)|\eta$, $\M^*=\M_\mu|(\Theta^+)^{\M_\mu}$ and whenever $\M^*\insegeq \N\insegeq \M_\mu$ is such that $\rho_\omega(\N)=\eta$, $\N$ is good.
\end{lemma}
\begin{proof}
We start by proving the first equality. We have that $\M^*\subseteq \H$ and therefore, $\M^*\in \H_\mu$. It is then enough to show that whenever $\M_\infty\insegeq \N\insegeq \M^*$ is such that $\rho_\omega(\N)=\Theta$, then $\N\in Ult(\H_\mu, \mu\cap \H_\mu)$. Since $\N$ is essentially a subset of $\Theta$ and $\mu\cap \H_\mu$ is a normal $\Theta$-complete ultrafilter in $\H_\mu$, we indeed have that $\N\in Ult(\H_\mu, \mu\cap \H_\mu)$. 


For the second equality, it is enough to show that whenever $\M_\infty\insegeq \N\insegeq \M_\mu$ is such that $\rho_\omega(\N)=\Theta$ then $\N$ is good. Let then $f:\Theta\rightarrow \M_\infty$ be such that $\pi_\mu(f)(\Theta)=\N$. It then follows that for a $\mu$-measure one set of $\k$, there is an elementary embedding $\pi_\k: f(\k)\rightarrow \N$. It then also follows that whenever $\pi:\bar{\N}\rightarrow \N$ is an elementary embedding such that $\bar{\N}$ is countable, there is a $\mu$-measure one set of $\k$ such that there is an embedding $\sigma: \bar{\N}\rightarrow f(\k)$ with the property that $\pi=\pi_\k\circ \sigma$. It then follows that $\bar{\N}$ is $\Sigma^\pi$-mouse that $\omega_1$-iterable. 

The proof of the third clause is just like the proof of the second.
\end{proof}

Next, we show that $\mu$ is amenable to $\M_\mu$.

\begin{lemma}\label{amenability}
$\mu$ is amenable to $\M_\mu$, that is $\mu \cap \N\in \M_\mu$ whenever $\N\inseg\M^*$.
\end{lemma}
\begin{proof}
This is a standard argument due to Kunen. Fix such an $\N$ projecting to $\Theta$. We have that $\pi_\mu\rest \N\in \pi_\mu(\M)$. The lemma now follows. 
\end{proof}

\subsection{A strategy for countable submodels of $L(\M_\mu, E_\mu)$}\label{realizable strategy sec}

Suppose $\xi$ is such that $E_\mu$ is a total extender in $L_\xi(\M_\mu, E_\mu)$. Suppose $\pi:\Q\rightarrow L_\xi(\M_\mu, E_\mu)$ is elementary and $\Q$ is countable. In this section, we show that $\Q$ has a $\pi$-realizable iteration strategy. Given any such $\Q$, we let $E^\Q$ be the preimage of $E_\mu$. We also let $\Q^-=\Q| ((\d^\Q)^+)^\Q$. Recall from \cite{ATHM} that if $\VT$ is a stack on some model $M$ and $R$ is a node in $\VT$ then $\VT_{\geq R}$ is the portion of $\VT$ after stage $R$. 

Fix then a $\Q$ as above. Given a stack $\VT$ on $\Q$, it can easily be partitioned into segments by considering when the image of $E^\Q$ is used. Thus, we  say $(\M_\a^0, E_\a, \M^1_\a, \VT_\a: \a<\eta)$ are the essential components of $\VT$ if 
\begin{enumerate}
\item $\M_0^0=\Q$,
\item $E_0$, $\M^1_0$ are defined if and only if the first extender used in $\VT$ is $E^\Q$ in which case $\M^1_0=Ult(\M^0_0, E^\Q)$.
\item $\VT_0$ is the largest initial segment of $\VT$ that is based on $(\M_1^1)^-$.
\item If $\a+1<\eta$ then $\M_{\a+1}^0$ is the last model of $\VT_\a$. Then $E_{\a+1}=E^{\M_{\a+1}^0}$ and $\M^1_\a=Ult(\M^0_{\a+1}, E_{\a+1})$. Again, $\VT_{\a+1}$ is the largest portion of $\VT_{\geq \M^1_{\a+1}}$ that is based on $\M_{\a+1}^1$. 
\item If $\a<\eta$ is limit then $\M_\a^0$ is the direct limit of $(\M_\b^0: \b<\a)$ under the iteration embeddings. The rest of the objects are defined as in the successor case. 
\end{enumerate}

Suppose now that $\VT$ is a stack on $\Q$ with essential components $(\M_\a^0, E_\a, \M^1_\a, \VT_\a: \a< \eta)$. Suppose that we also have embeddings $(\pi_\a^0, \pi_\a^1: \a< \eta)$ such that 
\begin{enumerate}
\item $\pi_0^0=\pi$.
\item For $\a<\eta$ and $i\in 2$, $\pi^i_\a:\M_\a^i\rightarrow L_\xi[\M, E_\mu]$.
\item For $\a<\b<\eta$ and $i, j\in 2$, $\pi^i_\a=\pi^j_\b\circ \pi^{\VT}_{\M_\a^i, \M_\b^j}$ and $\pi^1_\a=\pi^1_\a\circ \pi^\VT_{\M_\a^0, \M_{\a}^1}$.
\item For $\a<\eta$, $\VT_\a$ is according to $\pi^1_\a$-pullback of $\Sigma$. 
\end{enumerate}
We then say $\VT$ is $\pi$-\textit{realizable} if there is $\vec{\pi}=_{def}(\pi_\a^0, \pi_\a^1: \a< \eta)$ witnessing the above clauses. 

Suppose then $\VT=(\M_\a^0, E_\a, \M^1_\a, \VT_\a: \a\leq \eta)$ is such that $(\M_\a^0, E_\a, \M^1_\a, \VT_\a: \a< \eta)$ is $\pi$-realizable as witnessed by $(\pi_\a^0, \pi_\a^1: \a< \eta)$. We would like to define embeddings $(\pi^0_\eta, \pi^1_\eta)$ such that
\begin{enumerate}
\item For $i\in 2$, $\pi^i_\eta:\M_\eta^0\rightarrow L_\xi[\M, E_\mu]$.
\item  For $\a<\eta$ and $i, j\in 2$, $\pi^i_\a=\pi^j_\eta\circ \pi^{\VT}_{\M_\a^i, \M_\eta^j}$ and $\pi^1_\eta=\pi^1_\eta\circ \pi^\VT_{\M_\eta^0, \M_{\eta}^1}$.
\end{enumerate}
We then will have that $\VT$, without its last stack, is $\pi$-realizable as witness by $(\pi_\a^0, \pi_\a^1: \a\leq\eta)$.

Suppose first $\eta=\a+1$. Let $(\S, \Lambda)\in \mathcal{F}$ be such that 
\begin{center}
$\sup(\pi_{\a+1}^1[\d^{\M_{\a+1}^1}])=\d^{\M_\infty(\S, \Lambda)}$.
\end{center}
The above equality simply says that the direct limit of all $\Lambda$-iterates of $\S$ reaches the ordinal mentioned on the left side. 
Since $\M_{\a+1}^1$ is countable, we can also require that 
\begin{center}
$\pi_{\a+1}^1[\M^1_\a]\subseteq rng(\pi^\Lambda_{\S, \infty})$.
\end{center}
Let then $k=_{def} (\pi^{\Lambda}_{\S, \infty})^{-1}\circ \pi_{\a+1}^1:\M_{\a}^1\rightarrow \S$. We then have that $\VT_\a$ is according to $k$-pullback of $\Lambda$. It follows then that letting $\R$ be the last model of $k\VT_\a$ there is an embedding $l$ such that\\\\
(1) $l: \M_\eta^0\rightarrow \R$ and $\pi^{k\VT_\a}\circ k=l\circ \pi^{\VT_\a}$.\\\\
We then set $\pi^0_\eta=\pi^{\Lambda_\R}_{\R, \infty}\circ l$. Notice now that\\\\
(2) $\pi^\Lambda_{\S, \infty}=\pi^{\Lambda_\R}_{\R, \infty}\circ \pi^{k\VT_\a}$.\\\\
(1), (2) and our choice of $(\S, \Lambda)$ imply that\\\\
(3) $\pi^1_{\a}=\pi^0_\eta\circ \pi^{\VT}_{\M_\a^1, \M_\eta^0}$.\\\\
(3) then implies that $\pi^0_\eta$ is as desired. 

To define $\pi^1_\eta$ we use countable completeness of $\mu$. First let $\zeta=o((\M_\eta^0)^-)$. Thus, $\zeta$ is the successor of $\d^{\M_\eta^0}$ in $\M_\eta^0$. For each $a\in \zeta^{<\omega}$, let $A_a=\cap \{\pi^0_\eta(A): (a, A)\in E^{\M_\eta^0}\}$. Fix now a fiber $f$ for the set $\{ (\pi^0_\eta(a), A_a): a\in \zeta^{<\omega}\}$. We can now define 
\begin{center}
$\pi^{1}_\eta([a, g]_{E^{\M^0_\eta}})=\pi^0_\eta(g)(f(a))$.
\end{center}
It is a standard argument to show that $\pi^{1}_\eta$ is as desired. It is now easy to show that 

\begin{lemma}\label{realizable strategy} There is an $(\omega_1, \omega_1)$ iteration strategy $\Lambda$ for $\Q$ such that whenever $\VT$ is a stack according to $\Lambda$ then $\VT$ is $\pi$-realizable. 
\end{lemma}

\subsection{$L(\M_\mu, E_\mu)$ is a hod premouse}

In this section we show that no level of $L(\M_\mu, E_\mu)$ projects across or to $\Theta$. 

\begin{lemma}\label{not projecting}
$\N=L(\M_\mu, E_\mu)$ is a hod premouse such that for every $\a$, $\rho_\omega(\N|\a)>\Theta$.
\end{lemma}
\begin{proof}
We start with the following claim. Later we will show that in fact for every $\alpha$, $\rho_\omega(L_{\alpha}(\M_\mu,E_\mu))\neq \Theta$.\\

\textit{Claim 1.}
For every $\alpha$, $\rho_\omega(L_{\alpha}(\M_\mu,E_\mu)) \geq \Theta$

\begin{proof}
Suppose not and fix the least $\alpha$ such that $\rho_\omega(L_{\alpha}(\M_\mu,E)) < \Theta$. Now fix a $\beta$ such that $\beta <\lambda^{\M_\mu}$ and $\theta_{\beta} \leq \rho_\omega(L_{\alpha}(\M_\mu,E)) < \theta_{\beta+1}$. Fix a hod pair $(\Q^*,\Lambda^*)$ such that $\M_{\infty}(\Q^*,\Lambda^*)=\M_\mu(\beta+2)$ and an elementary hull $\sigma: H \to L_{\xi}(\P(\mathbb{R}))[\mu]$ such that:

\begin{enumerate}

\item
$\xi = \a^{+17}$,

\item
$\vert H \vert=\aleph_0$

\item
$(\Q^*,\Lambda^*)\in rng(\sigma)$,

\item
$L_{\alpha}(\M_\mu,E)\in rng(\sigma)$.

\end{enumerate}

Now let $\bar{\Q}=\sigma^{-1}(L_{\alpha}(\M_\mu,E)))$ and let $\gamma=\sigma^{-1}(\beta)$. Then by elementarity we have that $\rho_\omega(\bar{\Q})< \theta^{\bar{\Q}}_{\gamma+1}$. Let 
\begin{center}
$\Q=\mathcal{C}^{\bar{\Q}}_\omega(p, \bar{\Q}(\gamma+1))$
\end{center} where $p$ is the standard parameter of $\bar{\Q}$. Notice that $\Q$ is $\d_{\gamma+1}^\Q$-sound. Let $l:\Q\rightarrow \bar{\Q}$ be the core embedding. Let $\Lambda$ be a $\sigma\circ l$-realizable strategy of $\Q$ (see \rlem{realizable strategy}). 

It follows from the branch condensation of $\Lambda^*$ that\\\\
(1) $\Lambda_{\Q(\gamma+1)}=\Lambda^*_{\Q(\gamma+1)}$.\\\\
 Consider the pointclass generated by $(\Q,\Lambda)$:\[\Gamma(\Q,\Lambda)=\{A: \exists (\vec{\U},\R)\in B(\Q,\Lambda)(A <_W Code(\Lambda_{\R,\vec{\U}}))\}.\] Since $\Q$ is $\d^{\Q}_{\gamma+1}$-sound\\\\
 (2) $\Q$ is $\text{OD}_{\Lambda_{\Q(\gamma+1)}}$.\\\\
 This is because $\Q$ is the unique $\d_{\gamma+1}^\Q$-sound $\Lambda_{\Q(\gamma+1)}$-hod mouse generating the pointclass $\Gamma(\Q,\Lambda)$\footnote{Suppose $(\R, \Phi)$ is another $\Lambda_{\Q(\gamma+1)}$-hod pair such that $\Gamma(\Q, \Lambda)=\Gamma(\R, \Phi)$ and $\R$ is $\d_{\gamma+1}^\Q$-sound. We can compare $(\R, \Phi)$ with $(\Q, \Lambda)$ to a common pair $(\S, \Psi)$ and obtain embeddings $i: \R\rightarrow \S$ and $j:\Q\rightarrow \S$ such that $i\rest \d_{\gamma+1}^\Q=j\rest \d_{\gamma+1}^\Q=id$. Since both $\R$ and $\Q$ are $\d_{\gamma+1}^\Q$-sound, it follows that both $i$ and $j$ are the core embeddings. It then follows that $i=j$ implying that $\R=\Q$. We then get that $\Psi=\Phi^i$ and $\Lambda=\Phi^j$, implying that $\Psi=\Phi$.}. 
 
 Notice now that \\\\
 (3) $\bar{\Q}(\gamma+2)$ is $\Lambda_{\Q(\gamma+1)}$-full. \\\\
 This is simply because $\bar{\Q}(\gamma+2)$ is a $\Lambda^*$-iterate of $\Q^*$ and (as $\Q(\gamma+1)=\bar{\Q}(\gamma+1)$)
 \begin{center}
  $\Lambda_{\Q(\gamma+1)}=\Lambda_{\bar{\Q}(\gamma+1)}=\Lambda^*_{\Q(\gamma+1)}$. 
  \end{center}
  It now follows from (2) and (3) that we must have $\Q \in \bar{\Q}(\gamma+2)$. But $\rho_\omega(\Q)< \d^{\bar{\Q}}_{\gamma+1}$, contradiction!
\end{proof}

The following claim is easier and finishes the proof of the lemma.\\

\textit{Claim 2.}
For every $\alpha$, $\rho_\omega(L_{\alpha}(\M_\mu,E)))>\Theta$

\begin{proof}
Suppose not. Let $\alpha$ be the least such that $\rho(L_{\alpha}(\M_\mu,E)))=\Theta$. Then let $A \subseteq \Theta$ be definable over $L_{\alpha}(\M_\mu,E)$ such that $A \notin L_{\alpha}(\M_\mu,E)$. Notice that for every $\kappa<\Theta$, $A \cap \kappa \in \M_\mu$. 

Let $S$ be the set of $\k<\Theta$ such that $\theta_\k=\k$. Given $\k\in S$ we let $\M_{\kappa} \unlhd \M_\mu$ be the least such that $A \cap \theta_{\kappa} \in \M_{\kappa}$. 

Let $j:\M_\mu \to Ult(\M_\mu,\mu\cap \M_\mu)$ be the ultrapower embedding. Let $j((\M_{\kappa}:\kappa \in S))=(\N_{\kappa}:\kappa \in j(S))$. But then $A=j(A) \cap \Theta \in \N_{\Theta}$. Notice that $\N_\Theta$ is good (see \rlem{lower part characterization of m}) implying that $\N_\Theta \unlhd \M_\mu$. It follows that $A \in \M_\mu$, contradiction!
\end{proof}
\end{proof}

\subsection{$\powerset(\bR)$ is symmetrically generic over $L(\M_\mu, E_\mu)$}

We start by recalling Vopenka algebra. We work in $M=L(\powerset(\bR))$. Let $\varphi$ be a formula and $s\in \Theta^{<\omega}$. Given $t\in \Theta^{<\omega}$,  define \[\mathcal{A}^{s}_{\phi,t}=\{ a: dom(\vec{a})=dom(s), a(i)\subseteq s(i)^{\omega}, L(\powerset(\mathbb{R}))\models \phi[a,t]\}.\] We write
\begin{center}
 $(\phi, t)\equiv_s (\psi, v)$ if and only if $\mathcal{A}^{s}_{\phi,t}=\mathcal{A}^{s}_{\psi,v}$.
 \end{center}
 We then let $[\phi, t]_s$ be the $\equiv_s$-equivalence class of  $(\phi, t)$. 

Next define \[\mathbb{Q}=\{(s,[\varphi,t]): s\in \Theta^{<\omega}, \varphi \text{ is a formula }, t\in \Theta^{<\omega}\}\] The ordering on $\mathbb{Q}$ is defined as follows: \[(v,[\phi,s])\leq (u,[\psi, t]) \text{ if and only if } u \vartriangleleft v,\mathcal{A}^{u}_{\psi,t} \subseteq \mathcal{A}^{v}_{\varphi,s} \text{ and } \mathcal{A}^{v}_{\varphi,s} \restriction dom(u)\subseteq \mathcal{A}^{u}_{\psi,t}\]
where $\mathcal{A}^{v}_{\varphi,s} \restriction dom(u)=\{ a\rest dom(u) : a\in \mathcal{A}^{v}_{\varphi,s}\}$.\\

Let \[\mathbb{P}=\{(s,a): s \in\Theta^{<\omega}, dom(a)=dom(s)\  \text{and}\  \forall i \leq n, a(i)\subseteq s(i)^{\omega}\}.\]
We set $(s, a)\leq _\mathbb{P} (t, b)$ if and only if $t\insegeq s$ and $b\insegeq a$.

Suppose $g\subseteq \mathbb{P}$ is $\mathbb{P}$-generic. Let $f_g=\cup_{(s, a)\in g}s$. It follows from genericity that $f:\omega\rightarrow \Theta$ is a surjection. Also, for $i\in \omega$ let $h_g=\cup_{(s, a)\in g}a$. Again, it follows from genericity that $h:\omega\rightarrow \cup_{\a<\Theta}\a^{\omega}$ is a surjection.

The following lemma is standard and is an adaptation of arguments due to Vopenka.

\begin{lemma}[Vopenka's lemma]\label{vopenkas lemma}
Suppose $(f,h)$ is $\mathbb{P}$-generic. Let 
\begin{center}
$G=\{(f\rest n, [\varphi,v]): n<\omega, (f\rest n, [\varphi,v])\in \mathbb{Q} \wedge h\rest n\in \mathcal{A}^{f\rest n}_{\phi, v}\}$.
\end{center}
 Then $G$ is $\mathbb{Q}/\M_\mu$-generic. Therefore, it is $L(\M_\mu, E_\mu)$-generic for $\mathbb{Q}$.
\end{lemma}

\begin{proof} We first show that for $\xi<\Theta$, we can bound the ordinal parameters used to define $OD^M$ subsets of $\powerset(\xi)$.\\

\textit{Claim 1.} In $M$, there is a function $F:\Theta^{<\omega}\rightarrow \Theta$ such that whenever $t\in \Theta^{<\omega}$ and $a$ are such that $dom(a)=dom(t)$, $a(i)\subseteq t(i)^{\omega}$ and $a$ is $OD^M$ then $a$ is $OD$ in $L(\Gamma_{F(t)})$ where $\Gamma_\b=\{A\subseteq \bR: w(A)<\b\}$.\\\\ 
\begin{proof}
For $A\subseteq \bR^2$, we say $(\bR, A)$ codes an $OD$ structure if $(\bR, A)$ is a well-founded, extensional model of some fragment of $ZFC$ and its transitive collapse is $OD$. Notice that by a standard Skolem hull argument, in $M$,  if $A\subseteq \bR^2$ and $(\bR, A)$ codes an ordinal definable structure then for some $\b<\Theta$, $(\bR, A)$ codes an ordinal definable structure in $L(\Gamma_\b)$. Fix now $\a<\Theta$. Then we must have that
\begin{center}
$\sup\{ \b: \exists A\subseteq \bR^2 (A\in \Gamma_\a \wedge \b\ \text{is least such that}  L(\Gamma_\b)\models ``(\bR, A)$ codes an $OD$ structure")$\}<\Theta$.
\end{center}
This is because otherwise we will have a function $G:\Gamma_\a\rightarrow \Theta$ which is unbounded. As $\Theta$ is regular, this is impossible. 

Let then 
\begin{center}
$F(t)=\sup\{ \b+\omega :  \exists A \subseteq \bR^2 (A\in \Gamma_{\sup(rng(t))+\omega} \wedge L(\Gamma_\b)\models ``(\bR, A)$ codes an $OD$ structure"$)\}$.
\end{center}
\end{proof}

Let then $D\in \M_\mu$ be a dense subset of $\mathbb{Q}$. We want to see that $G\cap D\not =\emptyset$. 
Let $E$ be the set of $(s, a)$ such that $dom(s)=dom(a)$, $s\in \Theta^{<\omega}$, for every $i<dom(s)$, $a(i)\in s(i)^\omega$ and for some $[\phi, v]$ such that $(s, [\phi, v])\in D$, $a\in \mathcal{A}^s_{\phi, v}$. 

We claim that $E$ is dense in $\mathbb{P}$. Fix then $(s, a)\in \mathbb{P}$. \\

\textit{Claim 2.}
\begin{center}
$\cup\{\mathcal{A}^t_{\phi, v}\rest dom(s) : s\insegeq t ,(t, [\phi, v])\in D\}=\{b: dom(b)=dom(s)$, for all $i<dom(s)$, $b(i)\in s(i)^\omega\}$.
\end{center}
\begin{proof}
Notice that $D\not \in M$, which makes the claim non-trivial. Consider the set $A=\{ \mathcal{A}^t_{\phi, v}\rest dom(s): (t, [\phi, v])\in D\}$. Each member of $A$ is $OD^M$. We want to see that the density of $D$ implies that $(\cup A)^c=\emptyset$. 

To see this let $\N\insegeq \M_\mu$ be such that $\rho_\omega(\N)=\Theta$ and $D\in \N$. Fix $\zeta> \Theta$ and let $\pi: H\rightarrow L_\zeta(\powerset(\bR))[\mu]$ be such that 
\begin{enumerate}
\item $\N, F, D\in rng(\pi)$,
\item $\cp(\pi)=\k=\theta_\k$,
\item $s\in \k^{<\omega}$.
\end{enumerate}
Let $(\bar{\N}, \bar{F}, \bar{D})=\pi^{-1}(\N, F, D)$. We have that $\bar{F}=F\rest \theta_\k$ and $\bar{D}=\{ (t, ([\phi, v])^{L(\powerset(\bR))^H}): (t, [\phi, v])\in D\}$. Let $\bar{A}=\pi^{-1}(A)$. Then $\bar{\N}\insegeq \M_\infty=V_\Theta^{\H^M}$. We then have that $\bar{D}\in OD^M$. Thus, $\bar{A}$ is $OD^M$. 

Clearly, $\cup \bar{A}\subseteq \cup A$. Suppose then $B=_{def}(\cup \bar{A})^c\not=\emptyset$. We have that $B$ is OD$^M$.
It follows that $B$ is $OD^{L(\Gamma_{F(s)})}$. Because $\bar{D}$ is dense in $H$, we have that $B=\emptyset$.
\end{proof}

It follows from Claim 2 that there is $(t, [\phi, v])\in D$ such that $a\in \mathcal{A}^t_{\phi, v}\rest dom(s)$. Let then $b\in \mathcal{A}^t_{\phi, v}\rest dom(s)$ be such that $b\rest dom(a)=a$. Then $(t, b)\in E$ and $(t, b)\leq_{\mathbb{P}} (s, a)$.

Since we now have that $E$ is dense, we can fix $n<\omega$ such that $(f\rest n, h\rest n)\in E$. Let $[\phi, v]$ be such that $h\rest n\in \mathcal{A}^{f\rest n}_{\phi, v}$. We then have that $(f\rest n, [\phi, v])\in G\cap D$. We leave it to the reader to verify that $G$ is a filter.
\end{proof}

The following is the main lemma of this section.

\begin{lemma}\label{symmetric extension} $V$ is a symmetric extension of $L(\M_\mu, E_\mu)$. In fact, 
\begin{center}$V=L(\M_\mu, E_\mu)(\Theta^\omega)$. \end{center}
\end{lemma}
\begin{proof}
Let $(f, h)$ be $\mathbb{P}$-generic. Let $G$ be as in \rlem{vopenkas lemma}. Then $\Theta^\omega\in L(\M_\mu, E_\mu)[G]$. It is a consequence of $AD^+$ that every set of reals $A$ is $(OD_t)^M$ for some $t\in \Theta^\omega$. Therefore, we have that
\begin{center}
 $\powerset(\bR)\subseteq \H^M(\Theta^\omega)\subseteq L(\M_\mu, E_\mu)(\Theta^\omega)$.
 \end{center}
 The reader can find more on the above displayed formula by consulting Section 2 of \cite{CLSSSZ}. It then follows that  $L(\M_\mu, E_\mu)(\powerset(\bR))=L(\M_\mu, E_\mu)(\Theta^\omega)$. 
 
  Because every set $A\in \powerset(\Theta)\cap  L(\M_\mu, E_\mu)(\Theta^\omega)$ is added to $L(\M_\mu, E_\mu)$ by a small forcing (in fact by the Vopenka algebra at some $\theta_\a<\Theta$), we have that $E_\mu$ has a canonical extension $E_\mu^+$ to  $L(\M_\mu, E_\mu)(\Theta^\omega)$ (see Theorem 2.4 of \cite{CLSSSZ}). Let $\nu=\{ A: (\Theta, A)\in E_\mu^+\}$. We then have that 
  \begin{center}
  $L(\M_\mu, E_\mu)(\Theta^\omega)\models  ``\nu$ is a normal $\bR$-complete measure on $\Theta$".
  \end{center}
  It then follows that $V=L(\M_\mu, E_\mu)(\Theta^\omega)$\footnote{In fact, $\nu=\mu$ but we do not need this.}.
\end{proof}

%
%
%
%

\subsection{Computation of $\H$}

We can now easily conclude that the model $L(\M_\mu, E_\mu)$ is independent of $\mu$. The following corollary is a simple consequence of \rlem{symmetric extension} and the fact that small forcing doesn't create new measures.

\begin{corollary}
Suppose $\nu$ is a normal $\bR$-complete measure on $\Theta$. Then $E_\nu \in  L(\M_\mu, E_\mu)\in  L(\M_\mu, E_\mu)$.
\end{corollary}

\begin{corollary}\label{uniqueness of the model} Suppose $\mu$ and $\nu$ are two normal $\bR$-complete measures on $\Theta$. Then $L(\M_\mu, E_\mu)=L(\M_\mu, E_\nu)$.
\end{corollary}
\begin{proof}
We have that $E_\nu\in L(\M_\mu, E_\mu)$ and $E_\mu\in L(\M_\nu, E_\nu)$, implying the conclusion. 
\end{proof}

Borrowing \rlem{uniqueness of measures} from the next subsection, there is a unique normal $\bR$-complete measure $\mu$ over $\Theta$ minimizing $\lambda_\mu$. It then follows from \rcor{uniqueness of the model} that $L(\M_\mu, E_\mu)\subseteq \H$.  It follows from \rlem{symmetric extension} that $\H\subseteq L(\M_\mu, E_\mu)$. Putting these two results together, we can now show that
\begin{center}
$\text{HOD}=L(\M_\mu,E_\mu)$.
\end{center}

\subsection{The uniqueness of minimal measures}\label{uniqueness of measures sec}

In this section, we show that there is a unique measure minimizing $\lambda_\mu$. 

\begin{lemma}\label{uniqueness of measures} Suppose $\mu$ and $\nu$ are two normal $\bR$-complete measures over $\Theta$ such that $\lambda_\mu=\lambda_\nu$. Then $\mu=\nu$.
\end{lemma}

We spend most of this subsection proving \rlem{uniqueness of measures}. Suppose that $\mu$ and $\nu$ are as in the hypothesis of the theorem. Notice that it is enough to show that $E_\mu=E_\nu$. This is because given $E_\mu=E_\nu$, \rlem{symmetric extension} will imply that $\mu=\nu$. 

 Notice first that $\M_\mu=\M_\nu$. Let then $\M=_{def}\M_\mu$ and $\N=_{def}L(\M, E_\mu, E_\nu)$. Because $\l_\mu=\l_\nu$, we have that $\N$ is a bicephalus. 

Let now $\xi>\l_\mu$ be such that $\N|\xi\models ZFC-Powerset+``$there exists a largest cardinal" and let $\pi: \Q\rightarrow \N|\xi$ be a countable hull of $\N|\xi$. Using the construction introduced in \rsubsec{realizable strategy sec}, we can build an iteration strategy $\Lambda$ for $\Q$ which is $\pi$-realizable. Let $(E, F)=\pi^{-1}(E_\mu, E_\nu)$. The idea now is simple. If we succeed comparing $(\Q, \Lambda)$ with itself then we will get a contradiction as it will show that $E=F$. 

There are several issues with the above idea. The problem is that we do not a priori know that $\Lambda$ is fullness preserving. It need not be fullness preserving, $\Gamma$-fullness preserving would suffice. However, it is not clear how to define $\Gamma$. A similar issue arises in \cite{ATHM} where one needs to compare two hod pairs whose strategies are not fullness preserving. There the issue is taken care of as follows.

Suppose $(\P, \Phi)$ is a hod pair such that $\l^\P$ is a limit ordinal. First recall from \cite{ATHM} that
\begin{center}
$I(\P, \Phi)=\{ (\VT, \R): \VT$ is a stack on $\P$ according to $\Phi$ such that $\pi^\VT$ exists and $\R$ is the last model of $\VT\}$\\
$B(\P, \Phi)=\{ (\VT, \S): \exists (\VT, \R)\in I(\P, \Phi) ( \S\inseg_{hod}\R)\}$.
\end{center}
Also recall that $\S\inseg_{hod}\R$ means that $\S$ is a hod mouse initial segment of $\R$, i.e., there is $\a< \l^\R$ such that $\S=\R(\a)$. Following Section 3.2 of \cite{ATHM}, we let
\begin{center}
$\Gamma(\P, \Phi)=\{ A\subseteq \bR: \exists (\VT, \R)\in I(\P, \Phi) (A\leq_W Code(\Phi))\}$.
\end{center}

In the above formula, $\leq_W$ denotes Borel reducibility and $Code(\Phi)$ is the set of reals coding $\Phi$. It is shown in \cite{ATHM} (see Theorem 3.27 of \cite{ATHM}) that under some general conditions, $\Phi$ is $\Gamma(\P, \Phi)$ fullness preserving. The conditions needed for this fact are as follows.
\begin{enumerate}
\item $\Phi$ has hull condensation.
\item For all $(\VT, \R)\in B(\P, \Phi)$ there is a hod pair $(\S, \Psi)$ an an elementary embedding $\sigma: \R\rightarrow \S$ such that $\Phi_{\R, \VT}=(\sigma$ pullback of $\Psi)$ and $\Psi$ has branch condensation and is fullness preserving.
\end{enumerate}
Investigating the proof of Theorem 3.27 of \cite{ATHM}, it is not hard to notice that the hull condensation of $\Phi$ is needed to infer that\\\\
(1) whenever $(\VT, \R)\in I(\P, \Phi)$, $\Gamma(\P, \Phi)=\Gamma(\R, \Phi_{\R, \VT})$.\\\\
Clearly (1) is an easy consequence, modulo copy constructions, of hull condensation of $\Phi$. 

We now get back to our case and continue with the pair $(\Q, \Lambda)$. The first few definitions generalize immediately. Thus we let
\begin{center}
$I(\Q, \Lambda)=\{ (\VT, \R): \VT$ is a stack on $\Q$ according to $\Lambda$ such that $\pi^\VT$ exists and $\R$ is the last model of $\VT\}$\\
$B(\Q, \Lambda)=\{ (\VT, \S): \exists (\VT, \R)\in I(\Q, \Lambda) ( \S\inseg_{hod}\R)\}$\\
$\Gamma(\Q, \Lambda)=\{ A\subseteq \bR: \exists (\VT, \R)\in I(\Q, \Lambda) (A\leq_W Code(\Lambda))\}$.
\end{center}

Condition 2 above is also easily seen to be satisfied. This is because we have defined $\Lambda$ to be a $\pi$-realizable iteration strategy. Given then a $(\VT, \R)\in B(\Q, \Lambda)$, letting $\R^+$ be the last model of $\VT$, we have an embedding $\sigma: \R^+\rightarrow \N|\xi$ such that $\Lambda_{\R, \VT}$ is the $\sigma$-pullback of $\Sigma_{\M(\a)}$ where $\a$ is such that $\sigma(\R)=\M(\a)$. We can then find a hod pair $(\S, \Phi)$ such that $\M_\infty(\S, \Phi)=\M(\a)$ and an embedding $k:\R\rightarrow \S$ such that $\Lambda_{\R, \VT}=(k$-pullback of $\Phi)$. 

The only issue is that $\Lambda$ may not have hull condensation. Here is how to go around this problem. Suppose $(\VT, \R)\in I(\Q, \Lambda)$. Let then 
\begin{center}
$\Gamma_\VT=\Gamma(\R, \Lambda_{\R, \VT})$.
\end{center}
Notice that it follows that if $(\VT_0, \R_0)\in I(\Q, \Lambda)$ and $(\VT_1, \R_1)\in I(\R_0, \Lambda_{\R_0, \VT_0})$ then 
\begin{center}
$\Gamma_{\VT_0^\frown\VT_1}\leq_W \Gamma_{\VT_0}$
\end{center}
Let then $(\VT, \R)\in I(\Q, \Lambda)$ be such that $\Gamma_{\VT}$ is $\leq_W$-minimal. Thus, \\\\
(2) whenever $(\VU, \S)\in I(\R, \Lambda_{\R, \VT})$, $\Gamma_\VT=\Gamma_{\VT^\frown \VU}$\\\\
(2) now plays the role of (1). We can now repeat the proof of Theorem 3.27 and show that $\Lambda_{\R, \VT}$ is $\Gamma_{\VT}$-fullness preserving. We can then compare $(\R, \Lambda_{\R, \VT})$ with itself, a fact that implies $\pi^{\VT}(E)=\pi^{\VT}(F)$. Thus, $E=F$. 

\section{The existence of minimal models and some remarks}

We left open the question whether there are minimal models of $\Theta ms$. Here we show that they indeed exist. As was seen in \rlem{weakly min models}, there are weakly minimal models (assuming the existence of models of $\Theta ms$).

\begin{theorem}\label{existence of minimal models} Suppose $V$ is a weakly minimal model of $\Theta ms$. Then there is a minimal model of $\Theta ms$.
\end{theorem}
\begin{proof}
First recall \rrem{remark}. Let $\mu$ be a normal $\bR$-complete measure over $\Theta$. We then have that $V=L(\M_\mu, E_\mu)(\Theta^\omega)$. Suppose $W\subsetneq V$ is another weakly minimal model of $\Theta ms$ such that $W\subseteq V$. It follows that $\powerset(\bR)^W=\powerset(\bR)^V$. Working in $W$, let $\nu$ be a normal $\bR$-complete measure over $\Theta$. We then have that $W=L(\M_\nu^W, E^W_\nu)(\Theta^\omega)$. 

Notice now that $\M_\nu^W\insegeq \M_\mu$. Suppose first that $\M_\nu^W=\M_\mu$. It now follows from \rlem{symmetric extension} that $E_\nu^W\in L(\M_\mu, E_\mu)$ and that $V\models ``\nu$ is a normal $\bR$-complete measure over $\Theta$". It then follows from \rcor{uniqueness of the model} that $L(\M_\mu, E_\mu)=L(\M_\nu^W, E_\nu^W)$ implying that $V=W$.

It follows that we must have that $\M_\nu^W\inseg \M_\mu$. The above discussion then shows that if $W\subseteq V$ is a weakly minimal model of $\Theta ms$ minimizing the height of all possible $\M_\nu^W$ then $W$ itself is a minimal model of $\Theta ms$.
\end{proof}

We finish by asking the following question. The question is whether one can compute the $\H$ of the minimal model of $AD_{\bR}+``\Theta$ is a strong cardinal". First assume $AD^+$. We say that $\Theta$ is a strong cardinal if for every $\l\geq \Theta$ there is a set $E=\{ (a, A): a\in \l^{<\omega} \wedge A\in [\k]^{\card{a}}\}$ such that 
\begin{enumerate}
\item for every $a\in \l^{<\omega}$, $E_a$ is an $\bR$-complete ultrafilter over $\k^{\card{a}}$,
\item for $a\subseteq b\in \l^{<\omega}$, letting $\pi_{b, a}$ be the projection map\footnote{Recall that $\pi_{b, a}(t)=(t_{i_0}, t_{i_1},..., t_{\card{a}-1})$ where letting $b=(b_0, b_1, ..., b_{\card{b}-1})$, $(i_0, ..., i_{\card{a}-1})$ are chosen in such a way that $a=(b_{i_0}, b_{i_1},..., b_{\card{a}-1})$.}, for any $A\in \k^{\card{b}}$, $A\in E_b \iff \pi_{b, a}[A]\in E_a$, and
\item for any $((a_i, A_i): i\in \omega)\subseteq E$ there is $f:\cup_{i<\omega} a_i\rightarrow \k$ such that for every $i$, $f[a_i]\in A_i$, and 
\item for some $a\in [\l]^{<\omega}$, there is a function $f\in HOD_{E, \powerset(\bR)}$ such that $f:\Theta^{\card{a}}\rightarrow L_\Theta(\powerset(\bR))$ and $[a, f]_E=V_\l^{HOD_{E, \powerset(\bR)}}$. 
\end{enumerate}
We then let $``\Theta strong"$ stand for the theory $AD_{\bR}+``\Theta$ is a strong cardinal" and say that $M$ is \textit{a minimal model of $\Theta strong$} if $M$ is a transitive model of $ZF$ such that $Ord\cup \bR\subseteq M$, $M\models \Theta strong$ and for every transitive $N\subsetneq M$ such that $Ord\cup \bR\subseteq N$, $N\models \neg \Theta strong$.
We conjecture that minimal models of $\Theta strong$ exist.

\begin{conjecture} Assume $AD^+$ and that there is a largest Suslin cardinal which is a member of the Solovay sequence. Then there is a minimal model of $\Theta strong$. 
\end{conjecture}

The hypothesis of the conjecture is known as LSA. The reader can find more about it by consulting \cite{hodlsa}. There it is shown that the hypothesis is consistent relative to a Woodin cardinal which itself is a limit of Woodin cardinals (see Theorem 13.1 of \cite{hodlsa}).

\begin{question} Suppose $M$ is a minimal model of $\Theta strong$. What is the fine structural form of $\H^M$?
\end{question}

One may also consider models of $AD_\bR$ in which $\Theta$ is a strong cardinal and the model itself has a distinguished extender sequence above $\Theta$. Such models will perhaps have the form $\M=L^{\vec{E}}(\powerset(\bR))$ where $\vec{E}$ is a good extender sequence. We may allow that some of the extenders in $\M$ have critical point $\Theta$. If in addition we require that $\M$ is \textit{minimal} then it should be possible to study the fine structural properties of $\H^\M$. However, one has to be a bit more careful. While the authors haven't done much towards completing the project mentioned above, at the moment it is hard to perceive a proof from large cardinals that there is an $\M$ as above in which $\Theta$ is a strong cardinal and $\M\models ``$there is a Woodin cardinal which is a limit of Woodin cardinals". It seems that one would need to wait until the theory of hod mice evolves into that region. Nevertheless, it seems that studying models of LSA will lead towards confirmation of the following conjecture. 

\begin{conjecture} Assume LSA. Then there is a model of $AD_\bR$ of the form $\M=L^{\vec{E}}(\powerset(\bR))$ in which $\Theta$ is a strong cardinal and $\M\models ``$there is a Woodin cardinal". 
\end{conjecture}

Another question that authors find interesting is whether there are connections between $\Theta strong$ and super compactness of $\omega_1$. In \cite{TrangSuper}, Nam Trang has shown that $\Theta ms$ is equiconsistent with $AD_{\bR}+``\Theta$ is regular"+$``\omega_1$ is $\Theta$ supercompact". It can then be asked whether $\Theta strong$ is equiconsistent with any theory stipulating that $\omega_1$ has high degree of supercompactness. 

\begin{question} Assume $\Theta strong$. Is there a model of $AD_{\bR}+``\Theta$ is regular" +$``\omega_1$ is $\Theta^+$-supercompact"?
\end{question}

We suspect that the theory $AD_{\bR}+``\Theta$ is regular" +$``\omega_1$ is $\Theta^+$-supercompact" has a high consistency strength. A natural attempt to produce a model of $AD_{\bR}+``\Theta$ is regular" +$``\omega_1$ is $\Theta^+$-supercompact" is as follows. Start with a model of $AD^+$ and suppose $\a$ is such that $\theta_\a<\Theta$. Let $\l=(\theta_\a^+)^\H$ and let $\Gamma=\{ A\subseteq \bR: w(A)<\theta_\a\}$. Let $\mu$ be the $\omega_1$-supercompactness measure on $\l$. By a result of Woodin, it is unique. Consider then the model $L(\l^\omega, \Gamma, \mu)$. The only problem is that $\powerset(\bR)\cap L(\l^\omega, \Gamma, \mu)$ maybe bigger than $\Gamma$. In fact, it is an unpublished theorem of the second author (but see \cite{hodlsa}) that in the minimal model of $LSA$, this intersection is bigger than $\Gamma$.

\begin{theorem} Assume $AD^+$ and that for some $\a$, $\theta_\a<\Theta$. Let $\l=(\theta_\a^+)^\H$ and $\Gamma=\{ A\subseteq \bR: w(A)<\theta_\a\}$. Then if $L(\l^\omega, \Gamma)\cap \powerset(\bR)=\Gamma$ then there is $B\subseteq \bR$ such that $L(B, \mathbb{R})\models LSA$.
\end{theorem}

However, Woodin, in unpublished work, showed that the theory $AD_\bR+``\omega_1$ is supercompact" is consistent relative to the theory $ZFC+``$ there is a proper class of Woodin cardinals that are limit of Woodin cardinals". Based on intuitions coming from the theory of hod mice in the region of LSA we conjecture the following.

\begin{conjecture} The following theories are equiconsistent.
\begin{enumerate}
\item $ZF+AD_\bR+``\Theta$ is regular"+$``\omega_1$ is supercompact".
\item $ZFC+``$there is a Woodin cardinal which is a limit of Woodin cardinals"
\end{enumerate}
\end{conjecture}

\bibliographystyle{plain}
\bibliography{RachidGrigor1}
\end{document}